\documentclass{amsart}
\usepackage{a4wide,amssymb}

\usepackage[arrow,matrix,curve]{xy}\SilentMatrices
\def\xyma{\xymatrix@M.7em}

\def\Id{{\sf Id}}

\def\Gr { \sf Groups}
\def\Ab { \sf Ab}
\newcommand{\oo}{\omega}
\def\Ga{\Gamma}
\newcommand{\tp}{\otimes}
\newcommand{\Ker}{{\sf Ker}}
\newcommand{\cok}{{\sf Coker}}

\newcommand{\Ext}{{\sf Ext}}
\newcommand{\Hom}{{\sf Hom}}

\let\x\times
\let\d\partial

\def\Z{{\mathbb Z}}
\newtheorem{Pro}{Proposition}
\newtheorem{Le}[Pro]{Lemma}
\newtheorem{The}[Pro]{Theorem}
\newtheorem{Co}[Pro]{Corollary}
\theoremstyle{definition}
\newtheorem{De}[Pro]{Definition}
\theoremstyle{remark}

\begin{document}

\title
{Quadratic functors and one-connected two stage spaces}

\author
{H.-J. Baues}

\author{T. Pirashvili}

\maketitle

\section{Introduction} Let $T:\Gr\to \Gr$ be a functor. For
each simplicial group $G_*$ one obtains a new simplicial group $T(G_*)$ by
applying the functor $T$ on $G_*$. We are especially interested in functors of
the following type
$$T=(-)\tp Q,$$
where $Q$ is a square group and the tensor product is defined in \cite{square}.
If $G_*$ is an appropriate model of the loop space of the two dimensional
sphere, then the assignment $T\mapsto B(TG_*)$ yields the functor:
$$e:{\sf SG}\to {\sf CW}(2,3)$$
from the category of square groups to the homotopy category of one-connected
three types (i.~e. spaces $X$ with $\pi_iX=0$ for all $i\not =2,3$). One can ask
what sort of 3-types one gets in this way?

It is relatively easy to get a necessary condition for an object $X\in {\sf
CW}(2,3)$ to be of the form $B(TG_*)$, $T=(-)\tp Q$ for a square group $Q$.
Namely $X$ has to be flat, meaning that the corresponding $k$-invariant
$k:\Ga(\pi_2(X))\to \pi_3(X)$ factors through the kernel of $$\pi_2(X)\tp
\pi_2(X)\to \Lambda^2(\pi_2(X)).$$ We hope that the following conjecture is
true.

\

{\bf Conjecture}. For any flat  object   $X\in {\sf CW}(2,3)$ there exists a
square group $Q\in \sf SG$ and an  isomorphism $X\cong e(Q)$ in ${\sf
CW}(2,3)$.

\

Our main result claims that the conjecture is true provided $\pi_2X\in \bf A$,
where $\bf A$ is the smallest class of abelian groups which is closed under
arbitrary direct sums and contains i) all cyclic groups, ii) all abelian groups
$A$ such that 2 is invertible in $A$ and iii) all abelian groups $A$ such that
$\Ext(A,Sym^2A)=0$, where $Sym^2A$ is the second symmetric power of $A$. It is
clear that then $\bf A$ contains all finitely generated abelian groups as well
as all free and all divisible abelian groups.

We also consider the corresponding stable problem. Let $X_*$ be a simplicial
group, which is a model of the loop space of $S^n$ and let $Q$ be a square
group. By taking the $(n+2)$-th stage of the Postnikov tower of $BT(X_*)$
with $T=(-)\tp Q$,  one obtains the functor
$$\underline{e}_n:{\sf SG}\to {\sf CW}(n,n+1), \ \ n\geq 3.$$

We prove that if $X=\underline{e}_n(Q)$, then $\pi_{n+1}(X)$ is a vector space
over $\Z/2\Z$, and conversely if two annihilates $\pi_{n+1}(X)$, then there exists
a square group $Q$ such that $e_n(Q)\cong X$.

Our approach  is to use presquare groups. They are gadgets classifying
quadratic functors from the category of finite pointed sets to the category of
groups. If  $F$ is such a functor, we have
$$BF(S^1)\in {\sf CW}(2,3)$$
Thus one obtains a functor ${\sf PSG}\to {\sf CW}(2,3)$ from the category of
presquare groups to the category ${\sf CW}(2,3)$. Our result claims that an
object $X\in {\sf CW}(2,3)$  is isomorphic to one
of the form $BF(S^1)$, where $F$ corresponds to a presquare group, iff
$X$ is flat. Similarly, taking the $(n+1)$-th stage of the Postnikov tower
of $BF(S^n)$  one obtains the functor
$${\sf PSG}\to {\sf CW}(n,n+1), \ \ n\geq 3,$$
and we prove that any object of  ${\sf CW}(n,n+1)$ lies in the essential image
of this functor. To pass from presquare groups to square groups we then develop
appropriate obstruction theory.

The question about precise relationship between square groups and algebraic
models of 2-types was posed to the authors by M.~Jibladze.

\section{On certain quadratic functors}
In this section we consider few quadratic functors
defined on the category $\Ab$ of abelian groups.
Let $F:\Ab\to\Ab$ be a functor with $F(0)=0$.
Let us recall that the functor $F$ is {\it additive} or {\it linear} if the natural projection
$$ F(X \oplus Y) \rightarrow F(X) \oplus F(Y)$$
is an isomorphism. Furthermore, $F$ is {\it
quadratic} if the {\it second cross-effect} $$F(X \mid Y) = \Ker
(F(X\oplus Y)\to F(X)\oplus F(Y))$$ as a bifunctor is linear in $X$ and $Y$. In this
case one has a natural decomposition
$$F(X\oplus Y)\cong F(X)\oplus F(Y)\oplus F(X\mid Y).$$

 \subsection{Universal quadratic functor}
Let $A$ and $B$ be abelian groups.
 A  map $f:A\to B$ is
called {\it quadratic } if  the \emph{cross-effect} $$(a\mid
b)_f:=f(a+b)-f(a)-f(b)$$ is linear in $a$ and $b$. It follows then that
$f(0)=0$.
 It is well known that for
any abelian group $A$ there is a universal  quadratic function $p:A\to P(A)$,
meaning that for any quadratic map $f: A\to B$ there exists a unique
homomorphism $h: P(A)\to B$ such that $f=h\circ p$. In this way one obtains
the functor $A\mapsto P(A)$, which has the following alternative description.
Let $I(A)$ be the augmentation ideal of the group algebra of $A$. Then one has
the isomorphism
$$P(A)\cong I(A)/I(A)^2$$
induced by $p(a)\mapsto (a-1)({\sf mod} \ I(A)^2)$ (see \cite{approximation}).
 The
following fact is well known \cite{passi1}:
\begin{Le}
For any abelian group $A$ one has the following short exact sequence
\begin{equation}\label{sympea}
\xyma{0\ar[r] &Sym^2(A)\ar[r]^j& P(A)\ar[r]^q& A\ar[r] &0}
\end{equation}
where $Sym^2$ is the second symmetric power, the homomorphisms $j$ and $q$
are given by
$$j(ab)=(a\mid b)_p=p(a+b)-p(a)-p(b)$$
$$q(p(a))=a.$$
\end{Le}

 It follows that the functor $P$ commutes with
filtered colimits and one has the following natural isomorphism
$$P(A\oplus B)\cong P(A)\oplus P(B)\oplus (A\tp B).$$
Furthermore, one has  isomorphisms (\cite{passi})
 $$P(\Z)\cong \Z\oplus \Z,$$
$$P(\Z/2\Z)\cong \Z/4\Z,$$
$$P(\Z/2^n\Z)\cong \Z/2^{n+1}\Z\oplus  \Z/2^{n-1}\Z,  \ \ n>1$$
and for any odd prime $p$, 
one has an isomorphism
$$P(\Z/p^n\Z)\cong
\Z/p^n\Z \oplus \Z/p^n\Z$$

For an abelian group we let
\begin{equation}\label{theta}
\theta(A)\in \Ext(A,Sym^2(A)),\end{equation}
be the element corresponding to the exact sequence (\ref{sympea}).

\begin{Le}\label{thetasym} The class $\theta(A)$ is
represented by the canonical symmetric 2-cocycle $f^s$, given by
$$f^s(a,b)= ab\in Sym^2(A), \ \ a,b\in A.$$
\end{Le}

{\it Proof}. The homomorphism $q:P(A)\to A$ has a set-section $p:A\to P(A)$
and the cocycle corresponding to this section is exactly $f^s$.

The class $\theta$ is nontrivial in general. For example one has
$\theta(\Z/2^n\Z)\not =0$. However one has
\begin{Le}\label{ertimeoredi}
If $2$ is invertible in $A$, then $\theta(A)=0$.
\end{Le}

{\it Proof}. Let $g:A\to Sym^2(A)$ be the map given by $g(a)=a^2$. Then
$$(a\mid b)_g=2ab$$ which shows that the coboundary of $\frac{g}{2}$ is $f^s$

\subsection{A functor $\Psi$}For an abelian
group $A$ we let $\Psi(A)$ be the kernel of the natural projection $A\tp A\to
\Lambda^2(A)$ from the second tensor power to the second exterior power. Thus
by the very definition one has the following exact sequence \begin{equation}
\xymatrix{0\ar[r]& \ \Psi( A)\ar[r]& A\tp A\ar[r] &\Lambda^2 (A)\ar[r]&0}
\end{equation}
In this way one obtains  the functor $\Psi:\Ab\to\Ab$. The functor $\Psi$
commutes with filtered colimits and one has the following natural isomorphism
$$\Psi(A\oplus B)\cong \Psi(A)\oplus \Psi(B)\oplus (A\tp B).$$
Furthermore, one has isomorphisms $$\Psi(\Z)\cong \Z, \ \ \Psi(\Z/n\Z)\cong
\Z/n\Z.$$

\subsection{Whitehead $\Ga$-functor}\label{whiteheadgamma}

The functor $\Psi$ is closely related with Whitehead $\Ga$-functor, which is
defined as follows. Let $A$ and $B$ be abelian groups. A  quadratic
map $f:A\to B$ is
called {\it homogeneous} if $f(-a)=f(a)$.
It follows then that
$$(a\mid a)_f=-(a\mid -a)_f=f(a)+f(-a)=2f(a).$$ Based on this
identity a simple induction argument shows  that
$$f(na)=n^2f(a).$$
 It is well known \cite{hanrycertain} that for
any abelian group $A$ there is a universal homogeneous quadratic function
$\gamma:A\to \Ga(A)$, meaning that for any homogeneous quadratic map $f: A\to
B$ there exists a unique homomorphism $h:\Ga(A)\to B$, such that $f=h\circ
\gamma$. The functor $A\mapsto \Ga(A)$ is known as {\it the Whitehead's
quadratic functor}. It is well known \cite{hanrycertain} that the functor $\Ga$
commutes with filtered colimits and one has the following natural isomorphism
$$\Ga(A\oplus B)\cong \Ga(A)\oplus \Ga(B)\oplus (A\tp B).$$
Furthermore, one has  isomorphisms
 $$\Ga(\Z)\cong \Z,$$
$$\Ga(\Z/2^n\Z)\cong \Z/2^{n+1}\Z$$
and for any odd prime one has an isomorphism
$$\Ga(\Z/p^n\Z)\cong
\Z/p^n\Z.$$

It follows that if $a\in A$ is of order $n$, then $\gamma(a)\in \Ga(A)$ is of
order $n$ if $n$ is odd and  it is of order $2n$, provided $n$ is even.

We have a natural homomorphism $\tau:\Ga(A)\to A\tp A$, given by
$\tau(\gamma(a))=a\tp a$. It is clear that the image of $\tau$ lies in
$\Psi(A)$ and in this way one gets a natural homomorphism
$$\tau':\Ga(A)\to \Psi(A).$$
It is well known that $\tau'$  is an epimorphism, moreover it is an isomorphism
provided $A=\Z$, or $A=\Z/n\Z$ with odd $n$. To identify the kernel of this map
we need additional notations.

For each $n\geq 1$ and each abelian group $A$ we put $t_n(A)=\{a\in A\mid
2^na=0\}$. Multiplication by 2 yields the natural transformation
$t_{n+1}\to t_n$. We have also an inclusion $t_{n-1}\hookrightarrow t_n$. Here
and elsewhere we assume that $t_0=0$. Thus one obtains a natural transformation
$t_{n+1}\oplus t_{n-1}\to t_n$, whose cokernel is denoted by $\Phi_n$. It
follows that $\Phi_n:\Ab\to \Ab$ is a well-defined additive functor, which
commutes with filtered colimits and
$$\Phi_n(\Z)=0, \  \ \Phi_n(\Z/p^k\Z)=0$$
if $p$ is an odd prime. It is also clear that
$$\Phi_n(\Z/2^k\Z)=0, \ \ {\rm if } \ k\not =n$$
and $$\Phi_n(\Z/2^n\Z)=\Z/2\Z.$$
\begin{Le}
For each $n\geq 1$ there is  a well-defined homomorphism $\iota_n: \
\Phi_n(A)\to \Ga(A)$ given by $\iota_n(a)=2^n\gamma(a).$
\end{Le}

{\it Proof}. If $a,b\in t_n(A)$, then one has
$$2^n\gamma(a+b)=2^n\gamma(a)+2^n\gamma(b)+ 2^n(a\mid b)_
\gamma=2^n\gamma(a)+2^n\gamma(b)$$ Here we used the fact that $(a\mid
b)_\gamma$ is linear in $a$, and therefore $2^n(a\mid b)_ \gamma= (2^na\mid b)_
\gamma=0$. On the other hand if $a=2b$, then $2^n \gamma(a)=2^n
\gamma(2b)=2^{n+2}\gamma(b)=0$, because $2^{n+1}b=0$. Similarly, if $a\in
t_{n-1}$, then $2^{n-1}a=0$ and therefore $2^n\gamma(a)=0$. Thus $\iota_n$
is a well-defined homomorphism.

The collection $\iota_n$, $n\geq 1$, defines the natural transformation
$\iota:\Phi\to\Ga$, where $\Phi=\bigoplus _{n\geq 1}\Phi_n$.

\begin{Pro} For any abelian group $A$ the kernel of the natural map
$$\tau':\Ga(A)\to \Psi(A)$$ is isomorphic to $\Phi(A)$, thus
one has an exact sequence:
$$\xymatrix{0\ar[r]& \ \Phi( A)\ar[r]^{\iota}
&\Ga(A)\ar[r]^{\tau}& A\tp A\ar[r] &\Lambda^2 (A)\ar[r]&0}$$
\end{Pro}

{\it Proof}.  Let  us observe that $\Ga\to\tp^2$ induces a monomorphism on the
second cross-effect and therefore $\Phi':=\Ker(\Ga\to\tp^2)$ is an additive
functor. To show that $\iota$ yields an isomorphism $\Phi\to \Phi'$ it suffices
to evaluate on cyclic groups, because both functors in question are additive
and preserve filtered colimits. Since both functors vanish on $\Z$ and on
$\Z/n\Z$ with odd $n$, we have to consider only the case, when $A=\Z/2^n\Z$.
Since $\Ga(\Z/2^n\Z)$ is the cyclic group of order $2^{n+1}$ generated by
$\gamma(1)$ it follows that $\Phi'(\Z/2^n\Z)$ is the cyclic group of order two
generated by $2^n\gamma(1)$. On the other hand $\Phi_k(\Z/2^n\Z)=0$, provided
$n\not =k$ and $\Phi_n(\Z/2^n\Z)=\Z/2\Z$ and the result follows.

\begin{Co}\label{psidanamod2si} For any abelian group $A$ the natural
transformation $\Ga(A)\to \Z/2\Z\tp A$ induced by $\gamma(a)\mapsto \ a({\sf mod}\ 2A)$
factors trough $\Psi(A)$.
\end{Co}

{\it Proof}. It suffices to note that the composite $\Phi_n(A)\to A/2A$ is
induced by $a\mapsto 2^na=0$, $a\in t_n(A)$.

The functors $P$ and $\Ga$ are related via the natural transformation
$\nu:P\to \Ga$, which is given by $\nu(p(a))= \gamma(a)$. Since any homogeneous
quadratic function is quadratic it follows that this transformation is an
epimorphism. To identify the kernel, let us observe that the map
$f:A\to P(A)$ given by $f(a)=p(a)-p(-a)$ is linear. Indeed we have
$$(a\mid b)_f=(a\mid b)_p-(-a\mid -b)_p=0$$
because $(-\mid -)_p$ is bilinear.
\begin{Le}
One has the exact sequence
$$\xyma{0\ar[r]& _2A\ar[r]&
A\ar[r]^f&P(A)\ar[r]^\nu &\Ga(A)\ar[r]&0}$$
where $f(a)=p(a)-p(-a)$ and $_2A=\{a\in A\mid 2a=0\}.$
\end{Le}
{\it Proof}.
It is clear that the transformation $\nu:P\to \Ga$ yields an
isomorphism on the second cross-effects. Thus the kernel of $\nu$ is linear.
It is clear that $f$ yelds the transformation from the identity functor to
the kernel of $\nu:P\to \Ga$ and since both functors $\Id$ and $\Ker  (\nu)$
are linear and preserve filtered colimits, it suffices to
observe that the result is true for a cyclic
group $A$.

\section{Algebraic models of one-connected two stage spaces}
\subsection{One-connected two stage spaces}  For any $n\geq 2$, we
let ${\sf CW}(n,n+1)$ be the homotopy category of such pointed CW-complexes $X$
that $\pi_iX=0$ for all $i\not = n,n+1$. If $X$ is an object of ${\sf
CW}(n,n+1)$, then $\pi_i(\Sigma X)=0$ if $i<n+1$, thus the $(n+2)$-th stage of
the Postnikov tower of $\Sigma X$ belongs to ${\sf CW}(n+1,n+2)$. This yields
the functor
$$P_{n+2}\Sigma:{\sf CW}(n,n+1) \to {\sf CW}(n+1,n+2)$$
which is known to be an equivalence of categories, provided $n\geq 3$.

This category is closely related to the category $\Pi(n,n+1)$ of
$k$-invariants \cite{AH}, whose objects are triples $(\pi_n,\pi_{n+1},k)$,
 where $\pi_{n}$ and $\pi_{n+1}$  are abelian groups
  and $k:\Ga_n(\pi_n)\to \pi_{n+1}$ is a homomorphism. Here for an abelian
  group $A$ and a natural number $n\geq 2$, we let $\Ga_n(A)$ be $\Ga(A)$ if
  $n=2$ and $\Z/2\Z\tp A$ if $n\geq 3$. A morphism $f$
from $(\pi_n,\pi_{n+1},k)$ to $(\pi_n',\pi_{n+1}',k')$ is a pair
$(f_n,f_{n+1})$, where $f_n:\pi_n\to \pi_n'$ and  $f_{n+1}: \pi_{n+1}\to
\pi_{n+1}'$ are homomorphisms of abelian groups such that the diagram
$$\xymatrix{\Ga_n(\pi_n)\ar[r]^k\ar[d]_{\Gamma_n(f_n)}& \pi_{n+1}\ar[d]^{f_{n+1}}\\
\Ga_n(\pi_n')\ar[r]_{k'}& \pi_{n+1}'}$$ commutes. Taking the nontrivial
$k$-invariant yields the functor
$$\kappa:{\sf CW}(n, n+1)\to \Pi(n,n+1), \ \ n\geq 2$$ which fits in the
following linear extension of categories \cite{AH}, \cite{BW}  $$0\to D_n\to
{\sf CW}(n,n+1)\to \Pi(n,n+1)\to 0,$$ where $D_n$ is a bifunctor on $\Pi(n,n+1)$,
given by
$$D_n((\pi_n,\pi_{n+1},k),(\pi_n',\pi_{n+1}',k'))=\Ext(\pi_n,\pi_{n+1}').$$
In particular $\kappa:{\sf CW}(n,n+1)\to \Pi(n,n+1)$ yields a
bijection on isomorphism classes of objects, moreover
 $\kappa$ reflects isomorphisms and is surjective on morphisms.
\subsection{Braided and symmetric categorical groups}
We also need the following well-known algebraic models for ${\sf CW}(n,n+1)$,
$n\geq 2$ \cite{condu},\cite{braided}, \cite{2typesiterated}.
\begin{De}\label{BCG}
A  \emph{braided categorical group}, shortly BCG, consists of the following data
$$C=(\d:C_{ee}\to C_e,\ \ \{-,-\}:C_e\x C_e\to  C_{ee})$$  where $C_e$ and
$C_{ee}$ are groups and  $\d$ is a homomorphism, while $\{-,-\}$ is a map such
that the following equalities hold for $x,y,z\in C_e$ and $a,b\in C_{ee}$.
$$\d \{x,y\}=x^{-1}y^{-1}xy$$
$$\{\d a, \d b\}=a^{-1}b^{-1}ab$$
$$\{\d a,x\}\{x, \d a\}=1$$
$$\{x,yz\}=\{x,z\}\{x,y\}\{y^{-1}x^{-1}yx,z\}$$
$$\{xy,z\}=\{y^{-1}xy^{-1},y^{-1}zy\}\{y,z\}.$$
A braided categorical group is called \emph{symmetric
 categorical group}, shortly SCG, if
\begin{equation}\label{sym.ganmatrba}
 \{x,y\}\{y,x\}=1.\end{equation}
 \end{De}
It follows that $\Ker(\d)$ is an abelian group and ${\sf Im}(\d)$ is a normal
subgroup of $C_e$ and $\cok(\d)$ is an abelian group. One puts
$$\pi_0^C:=\cok(\d), \ \ \ \pi_1^C:=\Ker(\d).$$
The BCG's and SCG's form categories in an obvious way. A morphism of BCG's
(resp. SCG's) is called \emph{weak equivalence} if it induces an isomorphism on
$\pi_i$, $i=0,1$. Let ${\sf Ho(BCG)}$ (resp. ${\sf Ho(SCG)}$) denote the
localization of ${\sf Ho(BCG)}$ (resp. ${\sf Ho(SCG)}$) with respect to weak
equivalences.

 Let us note that BCG's are termed reduced 2-modules in
\cite{2typesiterated}, while SCG's  are termed stable 2-modules in
\cite{2typesiterated}. Thanks to \cite{condu} one knows that the category of
braided categorical groups is equivalent to the category of such simplicial
groups $G_*$, that $N_iG_*=0$, if $i\not =1,2$. Here $N_*G_*$ denotes the Moore
normalization of $G_*$. Similarly the category of symmetric categorical groups
is equivalent to the category of such simplicial groups $G_*$, that $N_iG_*=0$,
if $i\not =n,n+1$ for a fixed $n>1$. Therefore the classifying space functor
induces the functors $$b_2:{\sf BCG}\to {\sf CW}(2,3)$$ and
$$b_n:{\sf SCG}\to {\sf CW}(n,n+1), \ n\geq 3$$
such that $\pi_nb_n(C)=\pi_0^C$ and $\pi_{n+1}b_{n+1}(C)=\pi_1^C$, $n\geq 2$.

The inclusion functor ${\sf SCG}\subset {\sf BCG}$ has the left adjoint functor
$\lambda:{\sf BCG}\to {\sf SCG}$, which is obtained by
$$\lambda(C)=(\d:C_{ee}'\to C_e,\{-,-\}:C_e\x C_e\to C'_{ee}),$$
where $C_{ee}'$ is the  quotient of $C_{ee}$ by the relation
(\ref{sym.ganmatrba}).

 Then the  functor $\lambda$ makes the following
diagram commute:
$$\xymatrix{{\sf BCG}\ar[r]\ar[d]^\lambda& {\sf CW}(2,3)\ar[d]^{P_4\Sigma}\\
{\sf SCG}\ar[r]& {\sf CW}(3,4). }$$ According to
\cite{braided},\cite{2typesiterated} one has the equivalence of categories
$${\sf
Ho(BCG)}\cong {\sf CW}(2,3)$$ and  $${\sf Ho(SCG)}\cong {\sf CW}(n,n+1), \ \
n\geq 3.$$

\section{Presquare groups}

\subsection{Quadratic functors on pointed finite sets}
Let $\Ga$ be the category of finite pointed sets and let $\Gr$ be the
category of groups. We consider functors $F:\Ga\to \Gr$ with the property
$F([0])=0$. Here and elsewhere $[n]$ denotes the set $\{0,\cdots,n\}$, with
basepoint $0$. The functor $F$ is {\it linear} if the map
$$(Fr_1,Fr_2) : F(X \vee Y) \rightarrow F(X) \times F(Y)$$
is an isomorphism, where $X \vee Y$ is the sum in the category
$\Ga$ and $r_1 : X \vee Y \rightarrow X,\, r_2 : X \vee Y
\rightarrow Y$ are the retractions. Furthermore, $F$ is {\it
quadratic} if the {\it second cross-effect} $F(X \mid Y) = \Ker
(F(r_1),F(r_2))$ as a bifunctor is linear in $X$ and $Y$.

Let $H\Z:\Ga\to \Gr$ be the functor which assigns to a pointed set $S$ the free
abelian group generated by $S$ modulo the relation $*=0$, where $*$ is the
basepoint of $S$. For any abelian group $A$, we let $HA:\Ga\to \Gr$ be the
functor given by $HA(S)=A\tp H\Z(S)$. It is clear that $HA$ is a linear
functor. It is easy to prove that any linear functor $\Ga\to \Gr$ is isomorphic
to one of the form $HA$. Thus the assignment $A\mapsto HA$ is an equivalence
between the category of abelian groups and the category of linear functors
$\Ga\to \Gr$. The category of quadratic functors $\Ga\to \Gr$ has the following
description \cite{doldann}.

\begin{De}\label{gan.psg} A presquare
group, shortly a PSG, consists of the following diagram
$$M=( \xymatrix{M_e \times M_e\ar[r]^{\{-,-\}} & M_{ee}\ar[r]^{\sigma}&
M_{ee}\ar[r]^P &M_e}),$$ where $M_{ee}$ is an abelian group and $\sigma$ is a
homomorphism with $\sigma^2=\Id$. Moreover, $M_e$ is a group written
additively, $P$ is a homomorphism and $\{-,-\}$ is a bilinear map, that is
$\{x+y,z\}=\{x,z\}+\{y,z\}$ and $\{x,y+z\}=\{x,y\}+\{x,z\}$, for all $x,y,z\in
M_e$. One requires that
\begin{enumerate}
\item[(a)]

$P\sigma=P$,

\item[(b)]

$\sigma\{x,y\}+\{y,x\}=0$, $x,y\in M_e$,

\item[(c)] $P\{x,y\}=x+y-x-y$, $x,y\in M_e$,

\item[(d)]

$\{x,Pa\}=0$,  $x\in M_e$,  $a\in M_{ee}$.
\end{enumerate}

\end{De}
It follows from (b) that for any PSG $M$ one has $\{Pa,x\}=0$. It follows from
(c) and (d) that $Pa$ lies in the centrum of $M_e$. Thus $\cok (P)$ is
well-defined and by (c) it is an abelian group. It follows that $M_e$ is a
group of nilpotency degree 2. It follows from the condition (a) that $\sigma$
yields a well-defined involution on $\Ker(P)$.

If $M$ and $N$ are two PSG, then a morphism $f$ from $M$ to $N$ consists of
a pair of homomorphisms $f_e:M_e\to N_e,f_{ee}:M_{ee}\to N_{ee}$ such that
$f_{ee}$ commutes with involutions and the diagrams
$$\xymatrix{M_{ee}\ar[r]^P\ar[d]_{f_{ee}}& M_e\ar[d]_{f_e} && M_e\times
M_e\ar[r]^{\{-,-\}}\ar[d]_{f_e\times f_e}& M_{ee}\ar[d]_{f_{ee}}\\
N_{ee}\ar[r]^P& N_e && N_e\times N_e\ar[r]^{\{-,-\}}& N_{ee}}
$$
commute. We let $\sf PSG$ be the category of presquare groups.

If $M$ is a PSG and $S$ is a pointed set with basepoint $*$, we let $S\odot M$
be the group generated by the symbols $s\odot x$ and $[s,t]\odot a$ with
$s,t\in S$, $x\in M_e$, $a\in M_{ee}$ subject to the relations
\begin{enumerate}
\item[]

$[s,s]\odot a=s\odot P(a)$

\item[]

$*\odot x=0=[*,s]\odot a $

\item []
$ [s,t]\odot a= [t,s]\odot\sigma(a)$

\item []
$[s,t]\odot \{x,y\}=-t\odot x-s\odot y+t \odot x+s\odot y$

\end{enumerate}
where $s\odot x$ is linear in $x$ and where $[s,t]\odot a$ is
central and linear in $a$.

 A result similar  to \cite{square} shows that the functor $S\mapsto
 S\odot M$ is a quadratic functor on $\Ga$ and  in this way one gets the
 equivalence between  the category $\sf PSG$ of  presquare groups and
 the category of quadratic functors from $\Ga$ to $\Gr$.
Actually this is a very particular case of much more general results obtained in
\cite{doldann}.

\subsection{Homotopy and $k$-invariant of a presquare group} Let $M$ be a PSG.
 We set $$\pi_1^M:=\Ker(P:M_{ee}\to M_e) \ \ {\rm
and} \ \ \pi_0^M:=\cok(P:M_{ee}\to M_e).$$ The involution $\sigma$ equips
$\pi_1^M$ with an involution, which is still denoted by $\sigma$.

 For any $x\in M_e$ we
let $\bar{x}$ be the class of $x$ in $\pi_0^M$. It follows from the condition
(d) of the definition of PSG, that  $\{-,-\}$ factors through $\pi_0^M$ and
thanks to (b) it yields the homomorphism
$$\{-,-\}:\pi_0^M\tp\pi_0^M\to M_{ee}^-,$$ where
$$M_{ee}^-:=\{a\in M_{ee}\mid a+\sigma(a)=0\}.$$
 We also need the homomorphism
$\oo:\Lambda^2(\pi_0^M)\to M_e$ which is induced by the commutator map:
$$\oo(\bar{x}\wedge \bar{y})=x+y-x-y.$$
Consider the following diagram:
$$\xymatrix{0\ar[r]& \ \Phi( \pi_0)\ar[r]^{\iota}
&\Ga(\pi_0)\ar@{-->}[d]\ar[r]^{\tau}& \pi_0\tp \pi_0\ar[r] \ar[d]_{\{-,-\}}&\Lambda^2
(\pi_0)\ar[r]\ar[d]^{\omega}&0&\\
& 0\ar[r]& \pi_1\ar[r]&M_{ee}\ar[r]^P&M_e\ar[r]& \pi_0\ar[r]&0 }$$ where
$\pi_i=\pi_i^M$, $i=0,1$. The diagram is commutative thanks to the property (c)
of the definition of PSG. Since the columns are exact, we see that there is a
well-defined morphism \begin{equation}\label{kainv}k^M=k:\Ga(\pi_0)\to
\pi_1\end{equation} given by $k(\gamma(\bar{x}))=\{x,x\}$. A diagram-chase
shows that $k\circ \iota=0$. Furthermore, the condition (b) of the definition
of PSG shows that the image of $k$ lies in $\pi_1^-=\{b\in \pi_1\mid
b+\sigma(b)=0\} $.

\subsection{Stable homotopy and stable $k$-invariant of a presquare group}
We let ${\sf PSG}_s$ be the full subcategory consisting of objects $M$ such
that the involution on $M_{ee}$ is trivial, that is  $\sigma(a)=a$ for all
$a\in M_{ee}$. In this case the bracket $\{-,-\}:\pi_0\tp \pi_0\to M_{ee}$
factors through $\tilde{\Lambda }^2(\pi_0)\to M_{ee}$, where $\tilde{\Lambda
}^2(A)$ is the quotient of $A\tp A$ by the relation $a\tp b+b\tp a\sim 0$.

 The inclusion ${\sf PSG}_s\subset {\sf PSG}$ has the left
adjoint given by $M\mapsto \underline{M}$, where
$$\underline{M}=( \xymatrix{M_e \times M_e\ar[r]^{\{-,-\}} &
M_{ee}/(\Id-\sigma) \ar[r]^{\Id}& M_{ee}/(\Id-\sigma)\ar[r]^P &M_e}).$$ The
fact that $P$ is still well-defined
 follows from the property (a) of Definition \ref{gan.psg}. Moreover, the
 quotient map $M\to \underline{M}$ is a morphism in category $\sf PSG$.
 We now put
$$\underline{\pi}_i^M:=\pi_i^{\underline{M}}, \ \ i=0,1.$$
Thus $\underline{\pi_0}=\pi_0$, while $\underline{\pi}_1=
\Ker(M_{ee}/(\Id-\sigma)\to M_e).$ In other words $\underline{\pi_i}^M$,
$i=0,1$ is the $i$-th homology of the following chain complex
\begin{equation}\label{q} Q_*(M):=
(\cdots \buildrel  \Id+\sigma \over \longrightarrow
M_{ee}\buildrel \Id-\sigma \over \longrightarrow M_{ee}\buildrel  \Id +\sigma
\over \longrightarrow M_{ee}\buildrel \Id-\sigma  \over \longrightarrow M_{ee}
\buildrel P \over \longrightarrow M_e)\end{equation} We define the homomorphism
$$\underline{k}:\Z/2\Z\tp \pi_0^M\to \underline{\pi_1}^M$$ by
$$\underline{k}(\bar{x}):=\{x,x\}({\sf mod}(\Id-\sigma)).$$
The homomorphism $\underline{k}$ fits in the following commutative diagram with
exact rows:
$$\xymatrix{0\ar[r]& \Z/2\Z\tp \pi_0\ar[r]^\nu\ar[d]^{\underline{k}}
&\tilde{\Lambda} ^2(\pi_0)\ar[r]\ar[d]&\Lambda^2
(\pi_0)\ar[r]\ar[d]^{\omega}&0&\\
 0\ar[r]& \underline{\pi_1}\ar[r]&M_{ee}/(\Id-\sigma)\ar[r]^P&M_e\ar[r]& \pi_0\ar[r]&0 }$$
where the homomorphism $\nu$ is induced by $\bar{x}\mapsto \bar{x}\tp
\bar{x}$, while the homomorphism $\tilde{\Lambda} ^2(\pi_0)\to
M_{ee}/(\Id-\sigma)$ is induced by $\{-,-\}$. The fact that the last map is
well-defined can be checked as follows:
$$\{x,y\}+\{y,x\}=-\sigma\{y,x\}+\{y,x\}\in {\sf Im}(\Id-\sigma).$$

The commutative diagram
$$\xymatrix{0\ar[r]& \pi_1\ar@{-->}[d] \ar[r] & M_{ee}\ar[r]^P\ar[d]&M_2\ar[d]^{\Id}\\
 0\ar[r]& \underline{\pi_1}\ar[r]&M_{ee}/(\Id-\sigma)\ar[r]^P&M_e
}$$ shows that there is a natural epimorphism $\epsilon^M:\pi_1^M\to
\underline{\pi}_1^M$, which is an isomorphism provided $M\in {\sf PSG}_s$.
\subsection{Product of presquare groups} The category $\sf PSG$ possesses
all limits and colimits. In the sequel we need the following explicit
construction of the product in $\sf PSG$. Let $M$ and $N$ be two PSG. Then
$(M\x N)$ is the PSG given by
$$(M\x N)_e=M_e\x N_e,$$
$$(M\x N)_{ee}=M_{ee}\x N_{ee},$$
$$\sigma(a,c)=(\sigma (a),\sigma(c)),$$
$$P(a,c)=(Pa,Pc),$$
$$\{(x,u), (y,v)\}=(\{x,y\}, \{u,v\}),$$
where $a\in M_{ee},c\in N_{ee}$, $x,y\in M_e, u,v\in N_e$.

It is clear that the functors $\pi_i$, $i=0,1, \underline{\pi_1}$
(and the morphisms $k$, $\underline{k}$ as well ) preserve the product.

\subsection{Coproduct of presquare groups} In the sequel we need
also the following explicit
construction of the coproduct in $\sf PSG$. But first we recall few facts on
the category $\sf Nil$ of nilpotent groups of class two. The inclusion functor
$$\sf Nil\to \Gr$$ has the left adjoint functor, which is given by the
\emph{nilization} functor:
$$G\mapsto G^{nil}=G/[G,[G,G]]$$
Let $G_1$ and
$G_2$ be two objects in $\sf Nil$, then the copruduct $G_1\vee G_2$ in
$\sf Nil$ is
obtained by the nilization of the free product of the groups $G_1$ and $G_2$.
It is well known that one has the following short exact sequence
$$0\to G^{ab}_1\tp G^{ab}_2\to G_1\vee G_2\to G_1\x G_2\to 0,$$
where $G^{ab}=G/[G,G]$ is the \emph{abelization} of $G$. This shows that
any element of $G_1\vee G_2$ (all groups are written additively)
can be written as a sum of elements $a+b+w$,
where $a\in G_1$, $b\in G_2$ and $w$ is a sum of commutators of the form
$a_1+b_1-a_1-b_1$, $a_1\in G_1$ and $b_1\in G_2$.

\begin{Le}\label{constnil} Let
$$0\to A\to X\overset{\eta}\longrightarrow G\to 0$$
and
$$0\to B\to Y\overset{\xi}\longrightarrow H\to 0$$
be central extensions in $\sf Nil$ with abelian $G$ and $H$.
 Define the group $Z$ as the quotient
$(X\vee Y)/\sim$, where the equivalence relation $\sim$ is generated by
$$a+y\sim y+a$$
$$b+x\sim x+b$$
where $a\in A$, $b\in B$, $x\in X$ and $y\in Y$. Then  one has the
following central extension of groups
$$0\to A\oplus B\oplus (G\tp H)\overset{j}\longrightarrow
Z\to
G\x H\to 0.$$
Here the homomorphism $j$ is given by $j(a+b+g\tp h)=a+b+(x+y-x-y)$, where
$x\in X$ and $y\in Y$ satisfy $\eta(x)=g$ and $\xi(y)=h$.
\end{Le}

{\it Proof}. It follows from the definition of the group $Z$ that $j$ is
a well-defined homomorphism, whose image is a normal subgroup of $Z$. It is
also clear that $\cok(j)\cong G\x H$. It remains to show that $j$ is
a monomorphism. To this end let us recall that for any abelian group
$M$ which is considered as a trivial $(G\x H)$-module one has
$$H^2(G\x H,M)\cong H^2(G,M)\oplus H^2(Y,M)\oplus \Hom(G\tp H,M).$$
We now take $M=A\oplus B\oplus (G\tp H)$ and we let $cl(Z)\in H^2(G\x H,M)$ be
the element whose components in the above decompositions are
$i_{1*}(cl(X)),i_{2*}(cl(Y)),i_3$. Here $cl(X)\in H^2(G,A)$ and $cl(Y)\in
H^2(H,B)$ are elements defined by the given central extensions, while $i_1:A\to
M$, $i_e:B\to M$ and $i_3:G\tp H\to M$ are standard inclusions. The class
$cl(Z)$ defines a central extension:
$$0\to A\oplus B\oplus (G\tp H )\to Z_1\to G\x H\to 0$$
Since $G$ and $H$ are abelian groups, it follows that $Z_1\in \sf Nil$.
By our construction one has a commutative diagram
$$\xyma{0\ar[r] & A\ar[r]\ar[d]& X\ar[r]\ar[d]&G\ar[r]\ar[d]&0\\
0\ar[r] & M\ar[r]& Z_1\ar[r]&G\x H\ar[r]&0\\
0\ar[r] & B\ar[r]\ar[u]& Y\ar[r]\ar[u]&H\ar[r]\ar[u]&0}
$$
Thus we have a canonical morphism $X\vee Y\to Z_1$ and one easily shows that
it yields the homomorphism $Z\to Z_1$ which makes the following diagram commute
$$\xymatrix{& A\oplus B \oplus (G\tp H)
\ar[r]\ar[d]^{\Id}& Z\ar[r]\ar[d]&G\x H\ar[r]\ar[d]&0\\
0\ar[r] &  A\oplus B \oplus (G\tp H)
\ar[r]& Z_1\ar[r]&G\x H\ar[r]&0.}
$$
It follows that $j:A\oplus B \oplus (G\tp H)\to Z $ is a monomorphism and the
proof is finished.

Now we construct coproducts in $\sf PSG$. Let
$$M = (M_e\x M_e
 \overset {\{-,-\}_M} \longrightarrow  M_{ee}\overset{\sigma_M}\longrightarrow
 M_{ee} \overset {P_M}  \longrightarrow M_e)$$
$$N = (N_e\x N_e
 \overset {\{-,-\}_N} \longrightarrow  N_{ee}\overset{\sigma_N}\longrightarrow
 N_{ee} \overset {N_M}  \longrightarrow N_e)$$
be  presquare groups. Let us recall that $M_e,N_e\in \sf Nil$.
The coproduct
$$M \vee N = ((M \vee N)_e\x (M\vee N)_e \overset {\{-,-\}} \longrightarrow
(M\vee N)_{ee}\overset {\sigma} \longrightarrow
(M \vee N)_{ee} \overset P \longrightarrow (M \vee N)_e)$$
in the category $\sf PSG$ is given by
$$(M \vee N)_{ee} = M_{ee} \oplus N_{ee} \oplus \pi_0^M\tp \pi_0^N
\oplus \pi_0^N \otimes \pi_0^M$$
$$(M \vee N)_e = (M_e \vee N_e) / \sim $$
Here the equivalence relation is generated by
$$P_M(x) + c \sim c + P_M(x),$$
$$a + P_N (u) \sim P_N (u) + a,$$
for
$x \in M_{ee}, \, c \in N_e, \, u \in N_{ee},\, a \in M_e.$
Let
$\bar{x} \in \pi_0^M,\, \bar{u} \in  \pi_0^N$
be the elements in cokernels represented by $x$ and $u$ respectively.
The operators $\sigma$ and $P$ for
$M \vee N$
are defined by
$$\sigma(x+u+\bar{a_1}\tp \bar{c}_1+ \bar c_2 \otimes \bar a_2)=
\sigma _M(x)+\sigma_N(u)+\bar{a_2}\tp\bar{c_2}+ \bar c_1 \otimes \bar a_1$$
$$P(x + u + \bar a_1 \otimes \bar c_1 + \bar c_2 \otimes \bar a_2) = P_M(x) + P_N(u) + (a_1 + c_1 -a_1 -c_1) + (c_2 + a_2 -c_2-a_2)$$
From this definition is it clear that $ \pi_0^{M\vee N}\cong \pi_0^M\oplus
\pi_0^N $. Now the map $\{-,-\}:\pi_0^{M\vee N}\tp \pi_0^{M\vee N}\to (M\vee
N)_{ee}$ is given by
$$\{\bar{a}+\bar{c},\bar{a}_1+\bar{c}_2\}=
\{\bar{a},\bar{a}_1\}_M+\{c_1,c_2\}_N+a\tp c_1+c\tp a_1.$$

\begin{Le}\label{konamravlispierti}
For any $M,N\in \sf PSG$ one has the following isomorphisms $$ \pi_0^{M\vee
N}\cong \pi_0^M\oplus \pi_0^N, \ \  \underline{\pi}_1^{M\vee N}\cong
\underline{\pi}_1^{M}\oplus \underline{\pi}_1^N$$
$$\pi_1^{M\vee N}\cong \pi_1^{M}\oplus \pi_1^N\oplus (\pi_0^M\tp \pi_0^N)
$$
\end{Le}

{\it Proof}. We already had the first isomorphism. To get other isomorphisms,
one has to apply Lemma \ref{constnil} to central extensions
$$0\to {\sf Im}(P_M) \to M_e\to \pi_0^M\to 0$$
and
$$0\to {\sf Im}(P_N) \to N_e\to \pi_0^N\to 0$$
to conclude that ${\sf Im} P_{M\vee N}\cong  {\sf Im}(P_M) \oplus
 {\sf Im}(P_N) \oplus (\pi_0^M\tp \pi_0^N)$ which obviously implies the result.

\subsection{A pushforward construction}\label{push} Let $M$ be a PSG and
let $f:\pi_1^M\to A$ be a
homomorphism of abelian groups with involutions. We can form the pushout
diagram in the category of abelian groups with involutions
$$\xymatrix{\pi_1^M\ar[r]^i \ar[d]^f &M_{ee}\ar[d]\\
A\ar[r]&f_*(M_{ee}).}$$ It follows from the properties of the pushout
construction that we have the following commutative diagram with exact rows:
$$\xymatrix{0\ar[r]&\pi_1^M \ar[r]^i \ar[d]^f &M_{ee}\ar[d] \ar[r]^P&M_e\ar[d]^{\Id}\ar[r]&
\pi_0^M\ar[d]^{\Id}\ar[r]& 0\\
0\ar[r]& A\ar[r]&f_*(M_{ee})\ar[r]^{f_*(P)}&M_e\ar[r]& \pi_0^M\ar[r]& 0}
$$
It is clear that $f_*(M)$ is also a $PSG$, where $f_*(M)_e=M_e$ and
$f_*(M)_{ee}=f_*(M_{ee})$ and the map $M_e\x M_e\to f_*(M_{ee})$ is the
composite of the map $M_e\x M_e\to M_{ee}$ and the homomorphism $M_{ee}\to
f_*(M_{ee})$.  Furthermore one has
$$\pi_0^{f_*(M)}=\pi_0^M,\ \ \pi_1^{f_*(M)}=A$$ and $k^{f_*(M)}=f\circ k^M$.
\subsection{Presquare groups and the universal coefficient theorem}\label{psguct}
In this section we construct a collection of presquare groups using the
universal coefficient theorem in group cohomology. Let us recall that for any
abelian groups $A$ and $B$ there is  a natural short exact sequence
$$\xymatrix{0\ar[r]& \Ext(A,B)\ar[r]& H^2(A,B)\ar[r]^{\sf c} &\Hom(\Lambda ^2(A),B)\ar[r]& 0}$$
which has a splitting  natural in $B$. Here we used the well-known
 isomorphism $H_2(A)\cong \Lambda ^2(A)$. The homomorphism $\sf c$ is given by the
 commutator map: If
 $$0\to B\to G\to A\to 0$$
 is a central extension, corresponding to an element $x\in H^2(A,B)$, then
 ${\sf c}(x):\Lambda ^2(A)\to B$ is given by $(a,b)\mapsto u+v-u-v$. Here $u$
 and $v$ are liftings of $a$ and $b$ to the group $G$ which is written
 additively.

Of special interest is the case when $B=\Lambda ^2(A)$. We let ${\sf T}_A$
be the set of equivalence classes of central extensions
$$(N_A)=(\xymatrix{0\ar[r]& \Lambda^2(A)\ar[r]^\mu &N_A\ar[r]^\rho &A\ar[r]&0})$$
such that ${\sf c}(N_A)=\Id_{\Lambda^2(A)}$. The set ${\sf T}_A$ is nonempty
and the group $\Ext(A,\Lambda^2(A))$ acts transitively and freely on  ${\sf
T}_A$.

If $(N_A)\in  {\sf T}_A$, then one can define the presquare group $\omega(N_A)$
as follows. By definition we put
$$\omega(N_A)_{e}=N_A,$$
$$\omega(N_A)_{ee}=A\tp A,$$
$$P(a\tp b)=\mu(a\wedge b),$$
$$ \sigma(a\tp b)=-b\tp a,$$
$$\{x,y\}=\rho(x)\tp \rho(y),$$
where  $x,y\in N_A$ and $a,b\in A.$ One easily checks that $\omega(N_A)$ is a
PSG.

By our construction we have:
\begin{Le}\label{omegashomotopia}
For any abelian group $A$ and any $(N_A)\in {\sf T}_A$ one has isomorphisms
$$\pi_0^{\omega(N_A)}\cong A,$$
$$\pi_1^{\omega(N_A)}\cong \Psi(A)$$
under which $k^{\omega(N_A)}$ corresponds to the canonical
homomorphism $\tau':\Ga(A)\to \Psi(A)$ induced by $$\tau:\Ga(A)\to A\tp A, \ \ x\mapsto x\tp x.$$
Moreover, additionally one has
$$\underline{\pi}_1^{\omega(N_A)}\cong \Z/2\Z\tp A$$
and $\underline{k}^{\omega(N_A)}=\Id_{\Z/2\Z\tp A}.$
\end{Le}

We can apply the functor ${\sf PSG \to PSG}_s$, $M\mapsto \underline{M}$  to
the presquare group $\omega(N_A)$. Here $A$ is any abelian group and $N_A\in
{\sf T}_A$. In this way one obtains an object $\underline{\omega}(N_A)\in \sf
PSG$. By definition one has
$$\underline{\omega}(N_A)_e=N_e, \ \ \underline{\omega}(N_A)_{ee}=
\tilde{\Lambda}^2(A),$$  the structure homomorphisms are given by
$\sigma=\Id_{\tilde{\Lambda}^2}$, $P(a\tilde{\wedge} b)=\mu(a\wedge b),$ and
$\{x,y\}=\rho(x) \tilde{\wedge}  \rho(y),$ where  $x,y\in N_A$ and $a,b\in A$,
compare with the definition of $\omega(N_A)$. By our construction we have:
\begin{Le}\label{omegaqvedaxaziani}
For any abelian group $A$ and any $(N_A)\in {\sf T}_A$ one has isomorphisms
$$\pi_0^{\underline{\omega}(N_A)}\cong A,$$
$$\pi_1^{\underline{\omega}(N_A)}\cong A/2A$$
which identify $k^{\underline{\omega}(N_A)}$ with the canonical transformation $\Ga(A)\to
A/2A$ induced by $ \gamma(x)\mapsto x ({\sf mod} \ 2A).$
Moreover, additionally one has
$$\underline{\pi}_1^{\underline{\omega}(N_A)}\cong\Z/2\Z\tp A $$
and $\underline{k}^{\underline{\omega}(N_A)}=\Id_{\Z/2\Z\tp A}$.
\end{Le}

\subsection{Presquare groups and braided categorical groups}
Forgetting the involution one gets  the functor
$$\Upsilon:{\sf PSG}\to {\sf BCG}
$$
which is given by
$$(\xymatrix{M_e \times M_e\ar[r]^{\{-,-\}} & M_{ee}\ar[r]^{\sigma}&
M_{ee}\ar[r]^P &M_e})\mapsto ((P:M_{ee}\to M_e),\{-,-\}).$$
 The same functor
can be obtained in terms of functors on $\Ga$  as follows.

Let us recall that there is a standard way to prolong a functor $F:\Ga\to \Gr$
to a functor from the category of pointed simplicial sets to the category of
simplicial groups $ {\sf s.Sets}_*\to {\sf s}.\Gr$. First using direct limits
one can prolong $F$ to a functor from the category of pointed sets $ {\sf
Sets}_*\to \Gr$, then by degreewise action one obtains a functor from the
category of pointed simplicial sets to the category of simplicial groups.
 By abuse of notation we will still denote
this functor by $F$. In particular one can use this construction for the
functor $F=(-)\odot M$ for a PSG $M$. In this paper we are particularly
interested in the evaluation of $F=(-)\odot M$ on simplicial spheres and
especially on $S^1$, which is the simplicial model of the circle with  two
nondegenerate simplices. Let us recall that $S^1$ is $[n]$ in dimension $n$.
Moreover $s_i:[n]\to [n+1]$ is the unique monotone injection whose image does
not contain $i+1$, while $d_i:[n]\to [n-1]$ is given by $d_i(j)=j$ if $j<i$,
$d_i(i)=i$ if $i<n$, $d_n(n)=0$ and $d_i(j)=j-1$ if $j>i$.

 \begin{Le}Let $M$ be a PSG and $F=(-)\odot M:\Ga\to \Gr$.
 Then the Moore complex associated to $F(S^1)$ is isomorphic to the following
 complex
 $$\cdots \to 0\to M_{ee}\buildrel P \over \longrightarrow M_e\to 0.$$

 \end{Le}

 \noindent {\it Proof}. The fact that the Moore complex associated to
$F(S^1)$ vanishes in dimensions $>2$ is a particular case of Proposition 5.9 of
\cite{doldann} and the computations in dimensions 1 and 2 are trivial (compare
also with (2.6) of \cite{qucf}).

\

  Since the Moore complex of $F(S^1)$ is trivial in all dimensions except dimensions one
  and two, it follows from  \cite{condu} that it corresponds to a BCG.
  Thanks to Lemma  this
 particular BCG is nothing but $\Upsilon(M)$. In particular $BF(S^1)$ has only
 two nontrivial homotopy groups
$\pi_2B(S^1)\cong \pi^M_0$ and $\pi_3B(S^1)\cong \pi_1^M$ and the unique
nontrivial $k$-invariant is given by the map $k_M\in
\Hom(\Ga(\pi_0^M),\pi_1^M)$ constructed in  equation (\ref{kainv}). It follows
that we have the following commutative diagram of categories and functors:
$$\xymatrix{{\sf BCG}\ar[rd]&\ar[l]_\Upsilon
{\sf PSG}\ar[r]^{\kappa^*}\ar[d]^{ev(S^1)}&\Pi^*(2,3)\ar[d]^{forgetful}\\
&{\sf CW}(2,3)\ar[r]_\kappa&\Pi(2,3)}
$$
Here $ev(S^1):{\sf PSG}\to {\sf CW}(2,3)$ is the functor which is given by
$$M\mapsto BF(S^1), \ \ {\rm where} \ \ F=(-) \odot M,$$
while the category $\Pi^*(2,3)$ is defined as follows. An object of the
category $\Pi^*(2,3)$ is a triple $(\pi_2,\pi_3,k)$,
 where $\pi_2$ is an abelian group, $\pi_3$ is an abelian group with involution
 $\sigma$
  and $k:\Ga(\pi_2)\to \pi_3^-$ is a homomorphism, where as usual we put
$$\pi_3^-:=\{a\in \pi_3\mid a+\sigma(a)=0\}.$$
If $(\pi_2,\pi_3,k)$ and $(\pi_2',\pi_3',k')$ are objects of $\Pi^*(2,3)$, then
a morphism $f$ from $(\pi_2,\pi_3,k)$ to $(\pi_2',\pi_3',k')$ is a pair
$(f_2,f_3)$, where $f_2:\pi_2\to \pi_2'$ is a  homomorphism of abelian groups,
while $f_3: \pi_3\to \pi_3'$ is a homomorphism of abelian groups with
involutions such that the diagram
$$\xymatrix{\Ga(\pi_2)\ar[r]^k\ar[d]_{\Gamma(f_2)}& \pi_3^-\ar[d]^{f_3}\\
\Ga(\pi_2')\ar[r]_{k'}& \pi_3^{'-}}$$ commutes. The functor $$\kappa^*:{\sf
PSQ}\to \Pi^*(2,3)$$ is given by $$\kappa^*(M)=(\pi_0^M,\pi_1^M,k^M).$$ We have
also the forgetful functor $\Pi^*(2,3)\to \Pi(2,3)$ which forgets the
involution on $\pi_3$. This functor has the retraction given by the inclusion
$\Pi(2,3)\hookrightarrow \Pi^*(2,3)$. Under this inclusion $\pi_3$ is
considered as a group with involution, given by $\sigma(a)=-a$.

\subsection{Realization of one-connected 3-types via presquare groups} In this section we
characterize objects of the categories ${\sf CW}(2,3)$ and $\Pi^*(2,3)$ which
are isomorphic to objects of the form $F(S^1)$ or $\kappa^*(M)$, where $M\in
\sf PSG$ and $F=(-)\odot M$.

 An object of $\Pi^*(2,3)$ (resp. $\Pi(2,3)$) is called {\it flat} if
the composite
$$\xymatrix{\Phi(\pi_2)\ar[r]^{\iota}& \Ga(\pi_2)\ar[r]^k&\pi_3}$$
is zero, where the functor $\Phi$ and the natural transformation $\iota$ were
defined in Section \ref{whiteheadgamma}, in other words $k$ factors trough
$\Psi(\pi_0)$. An object $X\in {\sf CW}(2,3)$ is called {\it flat} provided
$\kappa(X)$ is flat.

\begin{The} \label{flatreal} {\rm i)} The values of $\kappa^*$
(and therefore of $\kappa$ as well) are flat.

{\rm ii)}   Let $(\pi_2,\pi_3,k)$  be a flat object of the category
$\Pi^*(2,3)$. Then there exist $M\in \sf PSG$ and an isomorphism
$\kappa^*(M) \cong (\pi_2,\pi_3,k)$ in
$\Pi^*(2,3).$

{\rm iii)} An object $X\in {\sf CW}(2,3)$ is isomorphic to an object of
the form $F(S^1)$, with quadratic $F:\Ga\to \sf Groups$ iff $X$ is flat.

\end{The}

{\it Proof}. Part iii) is an immediate consequence of i) and ii) and 
properties of linear extensions of categories \cite{BW}. The first statement
follows from the diagram chase based on the following commutative diagram:
$$\xymatrix{0\ar[r]& \ \Phi( \pi_0)\ar[r]^{\iota}
&\Ga(\pi_0)\ar[r]^{\tau}\ar[d]_k& \pi_0\tp \pi_0
\ar[d]^{\{-,-\}}\\
& 0\ar[r]& \pi_1\ar[r]&M_{ee} }$$ For the second part we prove that a
pushforward construction applied on $\omega(N_A)$ does the job. Here $N_A$ is
any element of ${\sf T}_A$, where $A=\pi_2$. Indeed, we already observed that
$$\kappa^*(\omega(N_A))=(A,\Psi(A),\tau')$$
where the involution on $\Psi(A)$ is given by $z\mapsto -z$, $z\in \Psi(A)$.
Let us now take a flat object $(\pi_2,\pi_3,k)$
of the category $\Pi^*(2,3)$. It follows that one has the commutative diagram
$$\xymatrix{\Ga(\pi_2)\ar[r]^{\tau'}\ar[dr]^k&\Psi(\pi_2)\ar[d]^{k'}\\
& \pi_3
}$$
Thus one can take the pushforward construction
$M=k'_*(\omega (N_A))$, $A=\pi_2$. Then one has
$\kappa^*(M) \cong (\pi_2,\pi_3,k)$ in
$\Pi^*(2,3)$ and we are done.

\subsection{Presquare groups and symmetric categorical groups}
It is clear that $\Upsilon(M)$ is a symmetric categorical group provided $M\in
{\sf PSG}_s$. On the other hand one can take the composite of functors
$\Upsilon:{\sf PSG}\to {\sf BCG} $ and $\lambda:\sf BCG \to SCG$ to get the
functor
$$\lambda\circ \Upsilon:{\sf PSG}\to {\sf SCG}
$$
It is clear that
$$\lambda(\Upsilon(M))= ((P:M_{ee}/(\Id-\sigma)\to M_e),\{-,-\}).$$
Thus one has the following commutative diagram
$$\xyma{{\sf PSG}_s\ar[r]^i \ar[d]_\Upsilon &{\sf PSG}\ar[r]^j\ar[d]_
\Upsilon \ar[dr]_{\lambda \circ \Upsilon}&{\sf PSG}_s\ar[d]^\Upsilon\\
{\sf SGC}\ar[r]^{i_1}& {\sf BCG}\ar[r]_\lambda &{\sf SCG}}
$$
where $i$ and $i_1$ are the inclusions, while $j(M)=\underline{M}$.

Let us fix a natural number $n\geq 2$ and let
$S^n$ be a simplicial model of the
$n$-dimensional sphere, which has only two nondegenerate simplices. For any
functor $F:\Ga\to \Gr$ one obtains the simplicial group $F(S^n)$ by applying
the functor $F$ on $S^n$. If $F$ is quadratic, then the Moore normalization of
$S^n$ is trivial  in dimensions $>2n$ and $<n$ and it is isomorphic to
$$\to \cdots 0 \to Q_n(M) \to \cdots \to Q_0(M)\to 0\cdots \to 0$$
where $F=(-)\odot M$ and $Q_*(M)$ is defined in  (\ref{q}).
As we see for $n\geq
2$ the space $BF(S^n)$ in general does not belong to ${\sf CW}(n+1,n+2)$.
However one can take the $(n+2)$-th stage of the Postnikov tower of
$BF(S^n)$, which is denoted by $e_n(M)$. It follows that one has the
following commutative diagram of categories and functors: $$\xymatrix{{\sf
SCG}\ar[d]&\ar[l]_{\lambda\circ 
\Upsilon}
{\sf PSG}\ar[d]^{\underline {\kappa}}\ar[dl]^{e_n}\\
{\sf CW}(n+1, n+2)\ar[r]_{\kappa}&\Pi(n,n+1)}
$$
where $\underline{\kappa}:{\sf PSG}\to \Pi(n+1,n+2)$ is given by $M\mapsto
(\pi_0^M,\underline{\pi_1}^M,\underline{k}^M)$.
\begin{The} Let $n\geq 2$. For any
 object $X$ of the category ${\sf CW}(n+1,n+2)$ there exists
an object
 $M\in{\sf PSG}_s$
and an isomorphism $e_n(M) \cong X$ in ${\sf CW}(n+1,n+2).$
\end{The}

{\it Proof}. Since the functor $\kappa:{\sf CW}(n+1,n+2)\to \Pi(n+1,n+2)$
induces bijection on isomorphism classes of objects and realizes all morphisms
in $\Pi(n+1,n+2)$ it suffices to prove that for any  object
$(\pi_{n+1},\pi_{n+2},k)$ of the category $\Pi(n+1,n+2)$ there exists an object
$M\in\sf PSG$ and an isomorphism $\underline{\kappa}(M) \cong
(\pi_{n+1},\pi_{n+2},k)$ in the category $\Pi(n+1,n+2)$. The proof of this
statement is quite similar to the proof of Theorem \ref{flatreal}. Let us
recall that for any abelian group $A$ and any element $N_A\in {\sf T}_A$ in
Section \ref{push} we constructed $\underline{\omega}(N_A)\in {\sf PSG}_s$ with
the property
$$\underline{\kappa}(\underline{\omega}(N_A))= (A,A/2A, \Id_{A/2A}).$$
Take now any object $(\pi_{n+1},\pi_{n+2},k)\in \Pi(n+1,n+2)$, where
$k:\pi_n/2\pi_n\to \pi_{n+1}$ is a homomorphism. One can take the 
pushforward construction $k_*(\underline{\omega}(N_A))$, $A=\pi_n$ to get
an object of expected kind.

\section{Square groups}
\subsection{Quadratic functors on the category of finitely generated free groups}
We now consider functors $F:{\sf Gr_f}\to \Gr$, where $\sf Gr_f$ is the
category of finitely generated free groups. For groups $G_1$ and $G_2$, we let
$G_1* G_2$ be the coproduct in $\Gr$. The functor $F:{\sf Gr_f}\to \Gr$ is {
\it linear} if the map
$$(Fr_1,Fr_2) : F(X * Y) \rightarrow F(X) \times F(Y)$$
is an isomorphism, where $r_1 : X * Y \rightarrow X,\, r_2 : X *Y \rightarrow
Y$ are the retractions. Moreover $F$ is {\it quadratic} if $F(X \mid Y) =
\Ker(Fr_1,Fr_2)$ as a bifunctor is linear in $X$ and $Y$. The main result of
\cite{square}  shows that the category of such quadratic functors ${\sf
Gr_f}\to \Gr$ is equivalent to the category of square groups. Here {\emph{a
square group} is a diagram
$$Q= (Q_e \overset H \longrightarrow Q_{ee} \overset P \longrightarrow Q_e)$$
where $Q_{ee}$ is an abelian group and  $Q_e$ is a group. Both groups are
written additively. Moreover $P$ is a homomorphism and $H$ is a map such that
the cross effect
$$(x \mid y)_H := H(x+y)-H(x) -H(y)$$
is linear in $x,y \in Q_e$. In addition the following properties are satisfied
$$(Pa\mid x)_H =0,$$
$$P(x \mid y)_H = x+y-x-y, $$
$$PHP (a) = P(a) + P(a) ,$$
where $x,y \in Q_{e}$ and $a,b\in Q_{ee}$. It follows from the first two
identities that $P$ maps to the center of $Q_e$. The second equation shows also
that $\cok(P)$ is abelian. Hence $Q_e$ is a group of nilpotency degree 2.  For
square groups one has the following additional formulas (see \cite{square}):
$$(x \mid Pa)_H =0,$$
$$H(x+y-x-y) =-(y\mid x)_H + (x\mid y)_H.$$

Now we relate the square groups with presquare groups.

\begin{Le}\label{gapreeba} Let $Q$ be a square group. Then $$\wp(Q)=
(Q_e,Q_{ee},\sigma=HP-\Id,(-,-)_H,P)$$ is a presquare group.
\end{Le}
\noindent {\it Proof}. The axioms (b) and (c) of the definition of PSG hold by
the definition of square group. Let us observe that, $HP$ is a homomorphism
thanks to the identity $(Px\mid a)_H=0$. Thus one
has$$\sigma^2=HPHP-2HP+\Id=H(2P)-2HP+ \Id=\Id,$$ which shows that $\sigma$ is
an involution. We have also
$$P\sigma=P(HP-\Id)=PHP-P=P$$
and
$$\sigma (x\mid y)_H+(y,x)_H=HP(x\mid y)_H-(x,y)_H+(y,x)_H=0.$$
Here we used the identity $P(x \mid y)_H =x+y-x-y$ and known expression for
$H(x+y-x-y)$.

We let $\sf SG$ be the category of square groups.
The presquare group $\wp(Q)$ is called the
underlying PSG of a square group $Q$. By abuse of
notations we write $\pi_0^Q,\pi_1^Q, \underline{\pi}_1^Q$
and $k^Q$ instead off
$\pi_0^{\wp(Q)},\pi_1^{\wp(Q)}, \underline{\pi}_1^{\wp (Q)}$ and $k^{\wp(Q)}$.

\

Let
$G$
be a group and let
$Q$
be a square group. We define the group
$G \otimes Q$
by the generators
$g \otimes x$
and
$[g,h]\otimes a$
with
$g,h\in G, x \in Q_e$
and
$a \in Q_{ee}$
subject to the relations
$$
(g+h) \otimes x  = g \otimes x + h \otimes x + [g, h] \otimes H(x)$$
$$[g,g] \otimes a  = g \otimes P (a)$$
where
$g \otimes x $
is linear in
$x$
and where
$[g,h] \otimes a $
is central and linear in each variable
$g,h$
and
$a$.
In this way one gets a bifunctor

$$ \otimes : {\sf Gr_f} \times  {\sf SG}
\rightarrow \sf Groups$$

One can prove (\cite{square}) that in addition the following identities hold:
$$[g,h]\tp a=[h,g]\tp \sigma(a), \ \ \sigma=HP-\Id$$
$$-h\tp x-g\tp y+h\tp x+g\tp y=[g,h]\tp (x\mid y)_H.$$
For any $Q\in \sf SG$ the functor $(-)\tp Q:\sf Gr_f\to Groups$ is quadratic
and any quadratic functor is isomorphic to  $(-)\tp Q:\sf Gr_f\to Groups$ with
appropriate $Q\in \sf SG$ \cite{square}.

\

In terms of quadratic functors the relation between $(-)\tp Q$ and $(-)\odot
\wp(Q)$ can be seen as  follows. For a pointed set $S$ we let ${\langle
S\rangle} $ be the free group generated by $S$ modulo the relation $*=0$, where
$*$ is the base point of $S$. Then one has a natural isomorphism
$$S\odot \wp(Q) \ \cong \ {\langle
S\rangle} \tp Q.$$ In other words the following diagram commutes
$$\xymatrix{\Ga\ar[r]^{\langle - \rangle} \ar[rd]_{\odot \wp(Q)}&\sf{Gr_f}\ar[d]^{\tp Q}\\
& \sf{Groups}}$$

\subsection{Product and coproduct of square groups}
In the sequel we need the following explicit
construction of the product and coproduct  in $\sf SG$.
Let $M$ and $N$ be two SG. Then
$(M\x N)$ is the SG given by
$$(M\x N)_e=M_e\x N_e,$$
$$(M\x N)_{ee}=M_{ee}\x N_{ee},$$
$$ H(x,y)=(H_M(x),H_N(y))$$
$$P(a,c)=(P_Ma,P_Nc).$$
Thus the functor $\wp:\sf SG\to PSG$ commutes with products. As our next
construction shows it commutes also with coproducts. By abuse of notation
we denote the underlying presquare groups of $M$ and $N$ still by $M$ and $N$.
We can consider the coproduct $M\vee N$ in $\sf PSG$. Define
$$H:(M\vee N)_{e}\to (M\vee N)_{ee}$$
by $$H(x+u+(a_1+c_1-a_1-c_1))=H_M(a)+H_N(c)+\bar{a}\tp \bar{c}_1-\bar{c}_1\tp
\bar{a}_1$$
One checks that in this way one really gets the coproduct in $\sf SG$ (see
 7.11 of \cite{square}).

\subsection{Lifting of PSG} We are going to answer the following question.
For a given $M\in \sf PSG$ under what conditions does there exists a square
group $Q$ such that $\wp(Q)\cong M$? If such $Q$ exists it is called a lifting
of $M$.

It is easy to see that not all PSG have liftings. Indeed, take $M_e=0$ and
$M_{ee}=\Z$. We  let $\sigma$ be the trivial involution on $M_{ee}$ and $P=0$,
$\{-,-\}=0$. Then one obtains a PSG. This particular PSG is not of the form
$\wp(Q)$, because if $P=0$ in a square group, then $\sigma=HP-\Id=-\Id$. This
show that unlike the linear functors not any quadratic functor $\Ga\to \Gr$
factors through $\sf Gr_f$.

As the following easy lemma shows if a PSG $M$ has a lifting  $Q\in \sf SG$
such a lifting in general is not unique. In fact the set of liftings is a
torsor on an appropriate group.
\begin{Le}\label{dut} a) Let $Q$ be a square group and let $\alpha:\pi_0^Q\to Q_{ee}$ be a
homomorphism. We set
$$Q^{\alpha}_e=Q_e, \ \ Q^{\alpha}_{ee}=Q_{ee}, \ \ P^{\alpha}=P, $$
and $$H^{\alpha}(x)=H(x)+{\alpha}(\bar{x})$$ where $x\in Q_e$ and $\bar{x}$
denotes the class of $x$ in $\pi_0^Q$. Then $Q^{\alpha}$ is a square group and
$\wp(Q)=\wp(Q^{\alpha})$. Conversely, if $Q$ and $P$ are two square groups with
$\wp(P)=\wp(Q)$ then first of all $P_e=Q_e$ and $P_{ee}=Q_{ee}$, furthermore
there exists a unique homomorphism $\alpha:\pi_0^Q\to Q_{ee}$ such that
$P=Q^{\alpha}$.

b) Let $Q,Q'\in \sf SG$ and let $f_e:Q_e\to Q_e'$ and $f_{ee}:Q_{ee}\to
Q_{ee}'$ be homomorphism of groups such that $f=(f_e,f_{ee})$ defines the
morphism $\wp(Q)\to \wp(Q')$ in the category $\sf PSG$. Then there exists a
unique homomorphism $\alpha(f):\pi_0^{Q}\to Q_{ee}$ such that
$$H'f_e(x)=f_{ee}H(x)+\alpha(f)(\bar{x}), \ \ x\in Q_e.$$
In other words $\alpha(f)=0$ iff $f$ is a morphism in $\sf SG$.
\end{Le}

 We now consider the problem under what conditions an object $M\in
\sf PSG$ is isomorphic to one of the form $\wp(Q)$. Of course if such $Q$ 
exists then $Q_e=M_e$ and $Q_{ee}=M_{ee}$. Moreover the map $P$ in $Q$ is the
same as in $M$. Thus the problem is under what conditions does there exist $H$ with
appropriate properties.

\subsection{The category ${\sf PSG}_0$}
We let ${\sf PSG}_0$ be the full subcategory of the category $\sf PSG$ which
consists of such $M$ that
$$\pi_1^-=\pi_1.$$
In other words, one requires that if $Pa=0$ for an element $a\in M_{ee}$, then
$\sigma(a)=-a$.

\begin{Le} Let $Q\in \sf SG$. Then $\wp(Q)\in {\sf PSG}_0$.
\end{Le}

{\it Proof}. Let us recall that in $\wp(Q)$ the involution $\sigma$ is given
by   $\sigma=HP-\Id$. Thus, if $Pa=0$, then $\sigma(a)=-a$.

\begin{Le}\label{pataraashi} Let $M\in {\sf PSG}_0$. Then there exists the unique homomorphism
$$h:{\sf Im}(P)\to M_{ee}$$
such that $hP(a)=a+\sigma(a)$.
\end{Le}

{\it Proof}. Uniqueness is clear, because each element from ${\sf Im}(P)$ can
be written as $P(a)$. To prove existence, we have to show that if $Pa=Pb$ then
$a+\sigma(a)=b+\sigma(b)$. If this holds, then $a=b+c$ with $Pc=0$. Thus
$$a+\sigma(a)=b+c+\sigma(b)+\sigma(c)=b+\sigma(b).$$

\begin{Le} Let $M\in {\sf PSG}_0$. Then for the diagram
$$A=(M_{ee}\buildrel P \over \longrightarrow {\sf Im}(P)\buildrel h \over
\longrightarrow M_{ee})$$ one has $PhP=2P$ and $hPh=2$. In other words $A$ is
a quadratic $\Z$-module in the sense of \cite{quadfun}.
\end{Le}

{\it Proof}. For $a\in M_{ee}$ one has $PhP(a)=P(a)+P\sigma(a)=2P(a)$. On the
other hand we have $hPhP=h(2P)=2hP$. Since $P:M_{ee}\to {\sf Im}(P)$ is an
epimorphism it follows that $hPh=2h$.

\subsection{A cohomological obstruction for lifting}
To each object $M\in \sf PSG$ one can associate two cohomological invariants.
The first one is the class
$$[M_e]\in H^2(\pi_0,{\sf Im}(P)), \ \ \pi_0=\pi_0^M,$$
which is associated to the central extension of groups:
$$0\to {\sf Im}(P)\to M_e\to \pi_0\to 0.$$
The second one is the class
$$[M_{ee}]\in H^2(\pi_0,M_{ee})$$
which is represented by the 2-cocycle $f\in {\sf Z}^2(\pi_0,M_{ee})$, where
$$f(\bar{x},\bar{y})=\{x,y\}.$$

\begin{De} Let $M \in {\sf PSG}_0$. Define the class
$$\vartheta(M)\in H^2(\pi_0,M_{ee})$$
by
$$\vartheta(M):=[M_{ee}]-h_*([M_e]).$$
Here $h_*:  H^2(\pi_0,{\sf Im}(P))\to  H^2(\pi_0,M_{ee})$ is induced from
the homomorphism $h$ defined in Lemma \ref{pataraashi}.
\end{De}

\begin{The}\label{aweva} If $Q \in {\sf SG}$, then $\vartheta(\wp(Q))=0$. Conversely if
$M \in {\sf PSG}_0$ is an object with $\vartheta(M)=0$, then there exists a
square group $Q$ and an isomorphism $\wp(Q)\cong M$.
\end{The}

{\it Proof}. Take $M\in {\sf PSG}_0$. Let us choose a set section $s:\pi_0\to
M_e$ of the quotient $M_e\to \pi_0$. One can assume that $s(0)=0$. For any
$x\in M_e$ one has $x-s(x)\in {\sf Im}(P)$.
 The class $[M_e]$
is represented by the 2-cocycle $\xi$, which is defined by
$$s(\bar{x})+s(\bar{y}) =\xi(\bar{x},\bar{y})+s(\bar{x}+\bar{y}).$$
If $M=\wp(Q)$, then the map $h:{\sf Im}(P)\to M_{ee}$ is the restriction of $H$
to ${\sf Im}(P)$. We set $$g=H\circ s:\pi_0\to M_{ee}.$$ One has
$$H(x)=H(x-s\bar{x}+s\bar{x})=h(x-s\bar{x})+g(\bar{x})+(x-s(x)\mod s(x))_H
=h(x-s\bar{x})+g(\bar{x})
$$
because
$x-s\bar{x}=P(a)$ for some $a\in M_{ee}$ and $\{P(a),s\bar{x}\}=0=
(P(a)\mid s\bar{x})_H$. It
follows that
$$H(x+y)=h(x+y-s(\bar{x}+\bar{y}))+g(\bar{x}+\bar{y}).$$
Since $y-s(\bar{y})$ lies in the center of $M_e$ one can write
$$h(x+y-s(\bar{x}+\bar{y}))= h(x+y-s(\bar{y})
-s(\bar{x})+\xi(\bar{x},\bar{y}))=$$
$$h(x-s(\bar{x}))+h(y-s(\bar{y}))+h(\xi(\bar{x},\bar{y})),$$
because $h$ is a homomorphism. Thus one obtains
$$(x\mid y)_H=H(x+y)-H(x)-H(y)=$$
$$h(x-s(\bar{x}))+h(y-s(\bar{y}))+h(\xi(\bar{x},\bar{y})
+ g(\bar{x}+\bar{y}) -h(x-s(\bar{x}))-g(\bar{x})-h(y-s(\bar{y}))- g(\bar{y})=$$
$$h(\xi(\bar{x},\bar{y})+
(\bar{x}\mid \bar{y}_g.$$ Since $(x\mid y)_H=(\bar{x}\mid \bar{y})_H$
represents the
class $[M_{ee}]$, and the function $(\bar{x}\mid \bar{y})_g$ is the coboundary of
$g$, we see that $\vartheta(M)=0$. Conversely assume $M\in {\sf PSG}_0$ is such
object that $\vartheta(M)=0$. The first condition defines the homomorphism
$h:{\sf Im}(P)\to M_{ee}$, while the second condition says that there exists a
function $g:\pi_0\to M_{ee}$ such that
$\{x,y\}=(\bar{x}\mid \bar{y})_g+h(\xi(\bar{x},\bar{y}).$ Now we can define
$H:M_e\to M_{ee}$ by $H(x)=h(x-s\bar{x})+g(\bar{x})$. One checks easily that
$M$ equipped with this $H$ is indeed a square group.

\subsection{Lifting problem for $\omega({N_A})$}

In this section we consider the problem whether for a given abelian group $A$
there exists an element $(N_A)\in {\sf T}_A$ such that $\omega(N_A)$ has a
lifting as a square group. The answer to this question depends entirely on
the element $\theta(A)\in \Ext(A,Sym^2(A)),$ which was defined in (\ref{theta}) via
the exact sequence (\ref{sympea}), or equivalently via
the canonical symmetric 2-cocycle $f^s$, given by
$f^s(a,b)= ab\in Sym^2(A)$,  $a,b\in A.$
Let us recall that $\theta(A)=0$
provided $2$ is invertible in $A$ ( see Lemma \ref{ertimeoredi}).

As an application of Theorem \ref{aweva}  we obtain the following

\begin{The}\label{omegasrealizacia} Let $A$ be an abelian group. If
there exists an element $(N_A)\in {\sf T}_A$ such that $\omega(N_A)\in \sf PSG$
has a lifting in $\sf SG$ then $\theta(A)=0$. Conversely, if $\theta(A)=0$,
then there exists an element $(N_A)\in {\sf T}_A$ such that $\omega(N_A)\in \sf
PSG$ has a lifting in $\sf SG$. In particular such a lifting exists provided
$2$ is invertible in $A$ or $\Ext(A, Sym^2A)=0$.
\end{The}

{\it Proof}. First of all let us observe that for any abelian group $A$ and any
$(N_A)\in {\sf T}_A$ one has $\omega(N_A)\in {\sf PSG}_0$. For $\omega(N_A)$ we
have ${\sf Im}(P)=\Lambda ^2(A)$ and the homomorphism $h:\Lambda ^2(A)\to A\tp
A$ is nothing but $h(a\wedge b)=a\tp b-b\tp a$. It follows then that the image
of $\vartheta(\omega(N_A))\in H^2(A,A\tp A)$ under ${\sf c}:H^2(A,A\tp A)\to
\Hom(A,A\tp A)$  is zero. Thanks to the  universal coefficient theorem one has
$\vartheta(\omega(N_A))\in \Ext(A,A\tp A)$. On the other hand one has also the
short exact sequence
$$0\to \Lambda^2(A)\to A\tp A\to Sym^2(A)\to 0$$
where the first arrow is $h$. Thus one has exact sequences
$$\Ext(A,\Lambda ^2(A))\to \Ext(A,A\tp A)\to \Ext(A,Sym^2(A))\to 0$$
and
$$H^2(A,\Lambda ^2(A))\to H^2(A,A\tp A)\to H^2(A,Sym^2(A))$$
Let us recall that $\vartheta(\omega(N_A))=[M_{ee}]-h_*([M_e])$. One observes
that the first term depends only on $A$ and does not depend on $(N_A)\in {\sf
T}_A$. It follows thus that the image of $\vartheta(\omega(N_A))\in H^2(A,A\tp
A)$ in $H^2(A,S^2(A)$ is the same as the image of $[M_{ee}]$ in $ H^2(A,A\tp
A)$. But $[M_{ee}]$ was represented by the cocycle $(a,b)\mapsto a\tp b$ and
therefore the image of $[M_{ee}]$ in $ H^2(A,A\tp A)$ lies in
$\Ext(A,Sym^2(A))$ and it coincides with $\theta(A)$. If $(N_A)\in {\sf T}_A$
is such element that the presquare group $\omega (N_A)$ has lifting, then
$\vartheta(\omega(N_A))=0$ and a fortiori $\theta(A)=0$. Conversely, assume
$\theta(A)=0$, then the exact sequence for ext groups shows that there is an
element $x\in \Ext(A,\Lambda^2(A))$ which maps to $\vartheta(\omega(N_A))$. But
$\Ext(A,\Lambda^2(A))$ acts on ${\sf T}_A$. Therefore  using $x$ we can correct
$N$ to obtain another element $N'\in {\sf T}_A$ such that
$\vartheta(\omega(N_A'))=0$ and we are done.

\begin{Co} If $\Ext(A, A\tp A)=0$ then the set
${\sf T}_A$  is a singleton and
$\omega(N_A)\in \sf PSG$ has a lifting in $\sf
SG$, where $(N_A)$ is the unique element of ${\sf T}_A$.
\end{Co}

{\it Proof}. Since $\Ext(A,-):\sf Ab\to Ab$ is right exact, it follows
that
$$\Ext(A,\Lambda^2(A))=0=\Ext(A,Sym^2(A)).$$
The first equation shows that ${\sf T}_A$ is a singleton, while the second
equations shows that such lifting exists.

\subsection{Lifting problem for $\underline{\omega}({N_A})$}

Let $A$ be an abelian group. We let $\underline{\theta}(A)\in
\Ext(A,Sym^2(A/2A))$ be the image of $\theta(A)\in \Ext(A,Sym^2(A))$ under the
canonical map $Sym^2(A)\to Sym^2(A/2A)$. The following is a straightforward
variation of the main result of the previous section:

\begin{Le} Let $A$ be an abelian group. If
there exists an element $(N_A)\in {\sf T}_A$ such that
$\underline{\omega}(N_A)\in \sf PSG$ has a lifting in $\sf SG$ then
$\underline{\theta}(A)=0$. Conversely, if $\underline{\theta}(A)=0$, then there
exists an element $(N_A)\in {\sf T}_A$ such that $\underline{\omega}(N_A)\in
\sf PSG$ has a lifting in $\sf SG$ .\end{Le}

{\it Proof}. The only difference is to use the exact sequence
$$\Lambda^2(A)\to \tilde{\Lambda} ^2(A)\to Sym^2(A/2A)\to 0$$
where the first map is induced by $a\wedge b\mapsto a\tilde{\wedge} b-
b\tilde{\wedge} a=2  a\tilde{\wedge} b$.

\subsection{Realization of one-connected 3-types via square groups}
We have the functors
$$\xymatrix{{\sf SG}\ar[r]^{\wp}& {\sf PSG}\ar[r]^\Upsilon& {\sf BCG}\ar[r]^
{b_2}& {\sf CW}(2,3).}$$ In this section we study the composite functor
$$e:{\sf SG\to CW}(2,3).$$ From the homotopy theoretic point of view  the
functor $e$ is the same as $Q\mapsto B((\Omega S^2)\tp Q)$. Here $\Omega S^2 $
is the simplicial group, which is obtained by the degreewise action of the
functor
$${\langle -\rangle} :\Ga\to \Gr$$ on $S^1$. Here $S^1$ is the standard simplicial model of
the circle with two nondegenerate simplices. The fact that this particular
simplicial functor is of the homotopy type of the loop space on the
two-dimensional sphere $\Omega S^2$ follows from the classical result of
Milnor. The fact that the functor $Q\mapsto B((\Omega S^2)\tp Q)$  is
isomorphic to the composite $b_2\circ \Upsilon\circ \wp$ follows from the
following isomorphism of simplicial groups:
$$(\Omega \ S^2)\tp Q\cong (S^1)\odot \wp(Q).$$
In this section we ask the following question: is every flat object of ${\sf
CW}(2,3)$ isomorphic to one of the form $B((\Omega \ S^2)\tp Q)$, where
$Q\in {\sf SG}$?

An object $(\pi_2,\pi_3,k)$ of  $\Pi(2,3)$ is called {\it realizable via SG} if
there exists a square group $Q$ and an isomorphism
$$\kappa(B(( \Omega \ S^2 )\tp Q))\cong (\pi_2,\pi_3,k).$$
 An object $X\in {\sf CW}(2,3)$ is called {\it realizable via SG} provided
$\kappa(X)$ is  realizable via SG. In other words $X$ is isomorphic to
$B((\Omega S^2)\tp Q)\cong ev(S^1)(\wp(Q))$, where $ev(S^1):\sf PSG\to CW(2,3)$
is the same as in Section 2.7.

\begin{Le} If $(\pi_2',\pi_3',k')$ and $(\pi_2'',\pi_3'',k'')$ are realizable via
SG, then $(\pi_2'\x \pi_2'',\pi_3'\x \pi_3'',k'\x k'')$  is also realizable via
SG.
\end{Le}

{\it Proof}. Indeed, if $Q'$ and $Q''$ realize $(\pi_2',\pi_3',k')$ and
$(\pi_2'',\pi_3'',k'')$ respectively, then $Q'\x Q''$ realizes $(\pi_2'\x
\pi_2'',\pi_3'\x \pi_3'',k'\x k'')$.

An abelian group $A$ is called {\it realizable via SG} provided
$(A,\Psi(A),\tau')$ is realizable via SG.

\begin{Le} If $\pi_2$ is  realizable via
SG, then any flat object of the form
 $(\pi_2,\pi_3, k)$ is also realizable via SG.
\end{Le}

{\it Proof}. Let $Q$ realize $(\pi_2,\Psi(\pi_2),\tau')$. Since $(\pi_2,\pi_3,
k)$ is flat, the homomorphism $k$ is the composite $k=k'\circ \tau'$, where
$k':\Psi(\pi_2)\to \pi_3$ is defined uniquely. Let us recall that in Section
\ref{push} we defined the pushforward construction for PSG's. It is clear that
 pushforward construction of a square group has a square group structure in an
 obvious way. It follows that $k'_*(Q)\in \sf SG$ realizes $(\pi_2,\pi_3, k)$ .

\begin{Le}\label{uh} Let $(A_i)_{i\in I}$ be a family of abelian groups. If each
 $A_i$ is  realizable via SG, then $\oplus _{i\in I}A_i$ is also realizable via SG.
\end{Le}

{\it Proof}. Assume  $Q_i$ realizes $A_i$. We claim that the coproduct of
$Q_i$ in the category of square groups realizes $\oplus _{i\in I}A_i$. Since
$\pi_i$ respects filtered colimits, it suffices to assume that $I$ is finite
and therefore without loss of generality one can assume that $I$ consists
of two elements. In this case the result follows from the isomorphisms
of Lemma \ref{konamravlispierti}.

\begin{Le}\label{witt} a) If $\theta(A)=0$, then $A$ is realizable via SG. In
particular any free abelian group, or any divisible abelian group  is
realizable via SG. Moreover, if 2 is invertible in an abelian group $A$, then
$A$ is realizable via SG.

b) For any $n\geq 1$ the group $\Z/2^n\Z$ is realizable via SG.

\end{Le}

{\it Proof}. If  $\theta(A)=0$  there exists $(N_A)$ such that $\omega(N_A)$
has a square group structure (see Theorem \ref{omegasrealizacia}) and this SG
realizes $A$. b) Let us consider the following square group:
$$Q_e=\Z^/2^{n+1}\Z=Q_{ee}$$
The homomorphism $P$ is multiplication by $2^n$. Define the quadratic
map $$H:\Z^/2^{n+1}\to \Z^/2^{n+1}$$ by $H(x)=x^2-x$. One easily checks that in
this way one obtains a SG which realizes $\Z^/2^{n}.$

Let $\bf A$ be the smallest class of abelian groups which is closed under
arbitrary direct sums and contains i) $\Z/2\Z$, ii) all abelian groups $A$ such
that 2 is invertible in $A$ and iii) all abelian groups $A$ such that
$\Ext(A,Sym^2A)=0$, where $Sym^2A$ is the second symmetric power of $A$. It is
clear that then $\bf A$ contains all cyclic groups, and hence all  finitely
generated abelian groups as well as all free and all divisible abelian groups.

\begin{Co} Let $X\in {\sf CW}(2,3)$ be a flat object. Then $X$ is realizable via SG
provided $\pi_2X\in \bf A$.

\end{Co}

\subsection{Realization of stable two-stage spaces via square groups}
 Now we consider the corresponding stable problem.
 Let us fix an integer $n\geq 3$. We let $\underline{e}_n$ be  the composite of the
 following functors:
$$\xymatrix{{\sf SG}\ar[r]^{\wp}& {\sf PSG}\ar[r]^-{e_n}& {\sf CW}(n,n+1).}$$
where $e_n$ is the composite of the following functors:
$$\xymatrix{{\sf PSG}\ar[r]^\Upsilon& {\sf BCG}\ar[r]^
\lambda&{\sf SCG}\ar[r]^-{b_n}& {\sf CW}(n,n+1).}$$ From the homotopy theoretic
point of view  the functor $\underline{e}_n$ is the same as $Q\mapsto
P_{n+2}B((\Omega S^{n})\tp Q)$. 
In this section we ask the following question:
What sort of objects of ${\sf CW}(n,n+1)$ are isomorphic to ones of the form
$\underline{e}_n(Q)$, where  $Q\in {\sf SG}$?

We start with few easy observations. We let ${\sf SG}_s$ be the full
subcategory of the category $\sf SG$ consisting of objects $Q$ such that
for any $a\in Q_{ee}$ one has $HP(a)=2a$.

\begin{Le} For any $Q\in {\sf SG}_s$ one has $\wp(Q)\in {\sf PSG}_s$.
\end{Le}

{\it Proof}. The involution $\sigma$ on $(\wp(M))_{ee}=M_{ee}$ is defined by
$\sigma=HP-\Id$. Thus $\sigma=\Id$ iff $M\in {\sf SG}_s$.

\begin{Le} Let $Q$ be a square group. Then there exists the unique square group
structure on $\underline{Q}$ such that the quotient map $Q\to \underline{Q}$ is
a morphism in $\sf SG$, where
$$(\underline{Q})_e=Q_e,\ \ (\underline{Q})_{ee}=Q_{ee}/(HP-2\Id)$$
Moreover the functor ${\sf SG\to SG}_s$ is the left adjoint functor to the
inclusion functor ${\sf SG}_s\subset {\sf SG}$
\end{Le}

{\it Proof} is immediate.

\begin{Co} For any $Q\in \sf SG$ the group $\underline{\pi}_1^Q$ is a vector space over
$\Z/2\Z$.
\end{Co}

{\it Proof}. By the definition we have
$\underline{\pi}_1^Q=\pi_1^{\underline{Q}}$. Since $\wp(\underline{Q})\in {\sf
PSG}_s$ it suffices to show that if $M\in {\sf PSG}_0\bigcap {\sf PSG}_s$ then
2 annihilates $\pi_1^M$. But by definition $M\in {\sf PSG}_0$ implies that the
involution on $\pi_1^M$ is multiplication by $(-1)$, while  $M\in
{\sf PSG}_s$ implies that the involution on $\pi_1^M$ is trivial, hence the
result.

We let ${\sf CW}(n,n+1)_*$ be the full subcategory of ${\sf CW}(n,n+1)$
consisting of spaces $X$ such that $\pi_{n+1}X$ is a vector space over
$\Z/2\Z$. Thus the values of the functor $\underline{e}_n$ lie in ${\sf
CW}(n,n+1)_*$.

\begin{The} For any object $X\in {\sf CW}(n,n+1)_*$, $n\geq 3$ there exists
an object $Q\in {\sf SG}_s$ and an isomorphism $\underline{e}_n(Q)\cong X$ in
${\sf CW}(n,n+1)$. Moreover, one can assume that $Q_e$ is an abelian group.
\end{The}
{\it Proof}. Still it suffices to realize objects like $(\pi_n,\pi_{n+1},k)$,
where $\pi_{n+1}$ is a vector space over $\Z/2\Z$. Using pushforward
construction it suffices to consider the universal case $(A,A/2A,k)$, where
$k:A\to A/2A$ is the canonical projection. We choose a basis $(b_i)_{i\in
I}$ of $A/2A$. Let $B$ a the free $\Z/4\Z$-module, with a basis
$(\tilde{b}_i)_{i\in I}$. We have canonical epimorphisms $\epsilon: A\to A/2A$
$\epsilon(a)=\bar{a}$ and $\varepsilon:B\to A/2A$,
$\varepsilon(\tilde{b}_i)=b_i$. It follows that one has the following exact
sequence
$$\xyma{0\ar[r]& A/2A\ar[r]^{\alpha}& B\ar[r]^\varepsilon &A\ar[r]& 0}$$
where $\alpha(b_i)=2\tilde{b}_i$. Let us consider the corresponding pullback
diagram
$$\xyma{C\ar[r]^\rho\ar[d]_\varrho& A\ar[d]^\epsilon\\
B\ar[r]_\varepsilon&A/2A}$$ It follows that one has the following exact
sequence
$$\xyma{0\ar[r]&A/2A\ar[r]^\iota&C\ar[r]^\rho&A\ar[r]&0.}$$
 We now put
$$Q_e=C,\ \ \ Q_{ee}=A/2A\oplus A/2A,$$
$$P(a,x)=\iota(x), \ \ a,x\in A/2A.$$
and $H(c)=(0,h(\varrho(c)))$, where $c\in C$ and $h:B\to A/2A$ is the quadratic
map uniquely defined by the conditions: $h(\tilde{b}_i)=0$ and
$(\tilde{b}_i\mid \tilde{b}_j)_h=0$, if $i\not = j$ and $(\tilde{b}_i\mid
\tilde{b}_i)_h=b_i.$ Here $i,j\in I$. A direct computation shows that in this
way one really gets a PSG which realizes $(A,A/2A,k)$.

\subsection{The transformation $\Delta$} The homotopy groups $\pi_i^Q$, $i=0,1$
and the stable homotopy group $\underline{\pi}_1^Q$ of a square group $Q$
depends only on the underlying presquare group $\wp(Q)$. In \cite{square}
a homomorphism $\Delta _Q:\pi_0^Q\to \pi_1^Q$ was constructed, which
defines the natural transformation of functors defined on $\sf SG$. Recall
that
$$\Delta(\bar{x})=HPH(x)+H(x+x)-4H(x), \ x\in Q_e.$$
Since $$\Delta P=HPHP+2HP-4HP=2HP-2HP=0$$
we see that $\Delta$ is well-defined.
Since $\sigma=HP-\Id$, one can rewrite
\begin{equation}\label{deltasigma}
\Delta(\bar{x})=\sigma(H(x))-H(x)+(x\mid x)_H
\end{equation}
Now it is clear that $\Delta$ is additive, because
$$(x\mid y)_{\Delta}=\sigma (x\mid y)_H-(x\mid y)_H+(x\mid y)_H+(y\mid
x)_H=0.$$
We have also
$$P\Delta=P\sigma H-PH+P(x\mid x)_H=0.$$
Thus $\Delta$ really defines the natural transformation $\pi_0\to \pi_1$.

It follows from the identity (\ref{deltasigma}) that the following diagram is
commutative:
$$\xymatrix{\pi_0\ar[r]^{\Delta}\ar[d]& \pi_1\ar[d]^{\epsilon}\\
\pi_0/2\pi_0\ar[r]_{\underline{k}}&\underline{\pi}_1
}$$

Let $\alpha:\pi_0^Q\to Q_{ee}$ be a homomorphism; according to Lemma \ref{dut}
we have also the square group $Q^{\alpha}$ which has the same underlying
presquare group as $Q$ and therefore the same homotopy groups as $Q$. One
easily sees that  
$$\Delta^{\alpha}=\Delta +\sigma \alpha-\alpha$$
which shows that $\Delta$ could not be constructed only in terms of  presquare
groups.

\begin{Le} Let $A$ be a finitely generated abelian group and let
$B$ be any abelian group. Furthermore let $f:A\to B$ be any homomorphism. Then
there exists a square group $Q$ such that $\pi_0^Q=A$, $\pi_1^Q=B$ and
$\Delta_Q=f$.
\end{Le}

{\it Proof}. Using pushforward construction it suffices to consider the
universal case $B=A$ and $f=\Id_A$. An abelian group $A$ is called
\emph{$\Delta$-realizable} if there exists a square group $Q$ such that
$\pi_0^Q=A=\pi_1^Q$ and $\Delta_Q=\Id_A$. Since $\pi_i:\sf SG\to Ab$, $i=0,1$
takes finite products to
finite products and $\Delta_{M\x N}=(\Delta_M,\Delta _N)$
 it suffices to show that any cyclic
group is $\Delta$-realizable. Assume $2$ is invertible in $A$. Then we have the
following square group
$$Q_e=A=Q_{ee}, \ \  P=0, \ \ {\rm and} \ \ \ H(a)=-\frac{a}{2}$$
which realizes $A$. The square group $\Z_{nil}$ realizes $\Z$, where
$$(\Z_{nil})_e=\Z=(\Z_{nil})_{ee}, \ \ P=0 \ \ {\rm and} \ \
H(a)=\frac{a^2-a}{2}$$
Finally the square group constructed in the proof of the part b) of Lemma
\ref{witt} realizes $\Z/2^n\Z$ for all $n\geq 1$.


\begin{thebibliography}{33}
\bibitem{AH} H.-J. Baues. Algebraic homotopy. Cambridge Studies in Advanced
Math. 54(1990), 84-92.

\bibitem{quadfun} H.-J. Baues. Quadratic functors and metastable homotopy. J.
Pure and Appl. Algebra. 94(1994), 49--107.

\bibitem{BW} H.-J. Baues and G.Wirsching, Cohomology of small categories.
J. Pure Appl. Algebra 38 (1985), 187--211.

\bibitem{2typesiterated} H.-J.Baues and D. Conduch\'e. On the 2-type of an
iterated loop space. Forum Math. 9(1997), 721--738.

\bibitem{BHP} H.-J. Baues, M. Hartl and T. Pirashvili.
 Quadratic categories and square rings. J. Pure Appl. Algebra 122
(1997), 1--40.

\bibitem{square} H.-J. Baues and T. Pirashvili.
Quadratic endofunctors of the category of groups. Adv. Math. 141 (1999),
167--206.

\bibitem{qucf} H.-J. Baues and T. Pirashvili.
A universal coefficient theorem for quadratic functors. J. Pure Appl. Algebra
148(2000),  1--15.

\bibitem{condu} D. Conduch\'e. Modules crois\'es g\'en\'eralis\'es de longueur
2. J. Pure Appl. Algebra. 34(1984), 155--178.

\bibitem{braided} A. Joyal and R. H. Street. Adv. Math. 102(1993), 20--78.


\bibitem{passi} I.B.S. Passi. Polynomial maps on groups.  J. Algebra  9  1968 121--151.

\bibitem{passi1} I.B.S. Passi. Polynomial functors.  Proc. Cambridge Philos. Soc.  66  1969 505--512.

\bibitem{approximation} T. Pirashvili. Polynomial approximation of ${\rm Ext}$ and ${\rm Tor}$ groups in functor categories.
Comm. Algebra 21 (1993), no. 5, 1705--1719.

\bibitem{doldann} T. Pirashvili.
Dold-Kan type theorem for $\Ga$-groups. Math. Ann. 318(2000), 277--298.

\bibitem{hanrycertain} J. H. C. Whitehead. A certain exact sequence.
Ann. Math. 52(1950), 52--110.


\end{thebibliography}
\end{document}